\documentclass[10pt,a4paper]{amsart}
\usepackage[utf8]{inputenc}
\usepackage[english]{babel}
\usepackage{amsmath}
\usepackage{amsfonts}
\usepackage{amssymb}
\usepackage{graphicx}
\usepackage{lmodern}
\usepackage[left=3.5cm,right=3.5cm,top=2.5cm,bottom=2.5cm]{geometry}
\usepackage[]{geometry}
\usepackage{tcolorbox}
\usepackage[ruled,vlined]{algorithm2e}
\include{pythonlisting}

\newcommand{\rhod}{\rho_\text{des}}
\DeclareMathOperator*{\argmin}{arg\,min}

\newcommand{\Z}{\mathbb Z}
\newcommand{\R}{\mathbb R}

\newcommand{\J}{\mathcal J}
\newcommand{\bF}{\bar F}
\newcommand{\Uad}{U_\text{ad}}
\newcommand{\T}{\mathbb{T}}

\newcommand{\calK}{\mathcal{K}}

\newcommand{\con}{u}
\newcommand{\h}{h}

\newcommand{\domain}{\Omega}

\newcommand*\di{\mathop{}\!\mathrm{d}}
\newcommand*\dd{\mathop{}\!\mathrm{d}}

\newtheorem{thrm}{Theorem}[section]
\newtheorem{lmm}{Lemma}[section]
\newtheorem{rmrk}{Remark}[section]
\newtheorem{ass}{Assumption}
\newtheorem{prblm}{Problem}

\title{Mean-field optimal control for biological pattern formation}
\author[M. Burger]{Martin Burger}
\author[L. M. Kreusser]{Lisa Maria Kreusser}
\author[C. Totzeck]{Claudia Totzeck}
\begin{document}
	\begin{abstract}
		We propose a mean-field optimal control problem for the parameter identification of a  given pattern. The cost functional is based on the Wasserstein distance between the probability measures of the modeled and the desired patterns. The first-order optimality conditions corresponding to the optimal control problem are derived using a Lagrangian approach on the mean-field level. Based on these conditions we propose a gradient descent method to identify relevant parameters such as angle of rotation  and force scaling which may be spatially inhomogeneous. We discretize the first-order optimality conditions in order to employ the algorithm on the particle level.  Moreover, we prove a rate for the convergence of the controls as the number of particles used for the discretization tends to infinity. Numerical results for the spatially homogeneous case demonstrate the feasibility of the approach.
	\end{abstract}	

	\maketitle
	
	\section{Introduction}
	In the past years, interacting particle or agent systems have been widely used to model collective behavior in biology, sociology and economics.  Among the many examples of  applications are biological phenomena such as animal herding or flocking \cite{BPTT,BPTTR20,DOCBC2006}, cell movement  \cite{GC2008}, as well as sociological and economical processes like opinion formation \cite{DMPW2009}, pedestrian flow dynamics \cite{BDFMW2014,CollisionA}, price formation \cite{BCMW2013}, robotics \cite{TFYE2015} and data science \cite{KW}.
	
	Most of these models start from interacting particle systems and encode `first principles' of biological, social or economical interactions, inspired by Newtonian physics with its classical dynamical systems of first or second order. A key feature  is the formation of global patterns even if the agents only interact with each other at a local scale. These patterns include consensus, polarisation and clustering. The qualitative results on pattern formation has been achieved by intensive studies of the limit for infinitely many particles.
	
	Interacting particle models have been extended to include control actions. The control of the dynamics has recently become an active research area \cite{ACFK2017,CBO2,DHL2017,FS2014,PRT2015,CBO1,CBOMTW}. The impact of the control on pattern formation has been studied both on the level of agents as well as in the mean-field limit, and has been applied successfully to a wide range of applications including traffic flow \cite{HKK07} and  herds of animals \cite{BPTT,BPTTR20}.
	
	While optimal control has mainly been studied for isotropic interacting particle models such as the  Cucker-Smale model \cite{PRT2015}, this paper focuses on control actions for a class of agent-based models with anisotropic interaction forces. 
	 For a large number $N$ of interacting cells with positions $x_j=x_j(t)\in\R^2$, ${j=1,\ldots,N},$ at time $t$ and interaction forces $ G$, we consider models of the form
	\begin{align}\label{eq:particlemodel}
	\frac{\di x_j}{\di t}=\frac{1}{N}\sum_{\substack{k=1\\k\neq j}}^N  G(x_j-x_k,T(x_j)),
	\end{align}
	equipped with initial data $x_j(0)=x_j^\text{in},~j=1,\ldots,N$,  for given $x_j^\text{in} \in\R^2,~j=1,\ldots,N$.  Here, $T$ denotes an underlying tensor field influencing the interaction force $F$ in addition to the distance vector $x_j-x_k$. This tensor field is responsible for anisotropic pattern formation and \eqref{eq:particlemodel} can be regarded as a prototype for understanding complex phenomena in nature.
	An example of this class of models is the K\"ucken-Champod model \cite{patternformationanisotropicmodel,Merkel} describing the formation of fingerprint patterns based on the interaction of  $N$  cells.  The interacting particle model  \eqref{eq:particlemodel} has been studied in \cite{patternformationanisotropicmodel,meanfieldfingerprints,stabilityanalysisanisotropicmodel,During2017}. In particular, the particles align in line patterns. For spatially homogeneous $T$, the stationary solution is given by straight line patterns, while for general $T$, more complex line patterns are observed in numerical simulations as stationary solutions.
	
	Fingerprint identification algorithms are of great importance in forensic science and are increasingly used in biometric applications. Due to data protection and privacy regulations, many of these algorithms are usually developed based on realistic synthetic fingerprint images which motivates the simulation of realistic patterns based on biological models. Using particle models like \eqref{eq:particlemodel} fingerprint patterns can be produced as stationary solution where the adjustment of one of the model parameters is related to the distances between the fingerprint lines. The motivation of this paper is to produce patterns with desired features, such as certain angles of rotation or scalings. We therefore propose to consider an optimal control problem constrained by \eqref{eq:particlemodel}. This allows us to estimate the model parameters for a given desirable stationary  pattern. In this framework, we focus on the estimation of the angle of rotation of the pattern and the strength of the interaction force at each point of the given domain. The spatially inhomogeneous control problem considered here is an inverse problem and can be regarded as the first step towards modeling complex fingerprint patterns with specific features in the future.
	
	The cost functional used for the parameter identification measures the Wasserstein distance of the  current states and the desired states and penalizes parameter settings that are too far from some given reference states. The algorithm proposed for the numerical simulations is based on projected gradient descent methods, where the gradient is computed using an adjoint approach. The first-order optimality conditions, required for evaluating the gradient, are derived with the help of a Lagrangian ansatz, similar to the one in \cite{BPTT}.
	Choosing the first-optimize-then-discretize approach pays off in the numerical implementation. In fact, we can use different time discretizations for the forward and the adjoint solver. This allows us to save a lot of effort, as the forward solver is implemented with an explicit Euler scheme and the linear and stiff adjoint system is solved implicitly.
	
	The main novelties of our approach are the modeling of a spatially inhomogeneous mean-field optimal control problem with periodic boundary conditions, and the treatment of the Wasserstein cost. The latter is, in particular, challenging from the numerical point of view. We show some simulation results on the particle level to demonstrate the feasibility of our approach and show a rate for the convergence of the controls as the number of cells tends to infinity.
	
	This paper is organised as follows. In Section \ref{sec:model}, we introduce the mean-field model and the considered spatially inhomogeneous optimal control problem in the macroscopic setting. The first-order optimality conditions are derived in Section~\ref{sec:foc}. In Section \ref{sec:focdiscrete}, we derive the discrete optimal control problem and the first-order optimality conditions by discretizing the forward and adjoint models. We also show the convergence of the discrete optimal control problem. The derived discretizations of the problem form the basis for our numerical schemes and we describe the resulting algorithm in Section \ref{sec:numericalschemes}. Numerical simulation  results for the spatially homogeneous case are shown in Section \ref{sec:numericalresults} before we conclude the paper in Section~\ref{sec:conclusion}.

 		\section{Description of the problem}\label{sec:model}
 		In this section, we introduce an interacting particle model on the torus that includes control actions. Starting from \eqref{eq:particlemodel}, we introduce the considered interaction forces first, then we pass to the continuum model and formulate the optimal control problem.

 		\subsection{Interaction forces}
 		For formulating an optimal control problem, we introduce spatially inhomogeneous interaction forces depending on a control variable $u$. We introduce a Hilbert space $\mathcal H(\R^2)$ of real-valued functions on $\R^2$ and  require that it is continuously embedded in $L^\infty(\R^2)$ and $C^0(\R^2)$. An example for $\mathcal H(\R^2)$ is given by the Sobolev space $H^{2}(\R^2)$ which is continuously embedded in $L^\infty(\R^2)$ and $C^0(\R^2)$ in two dimensions.
 		We define the space of controls as $U = \mathcal H(\R^2)\times \mathcal H(\R^2)$ and consider controls $u=(\theta,\eta)\in U$ with $\theta \in \mathcal H(\R^2)$ and $\eta \in \mathcal H(\R^2)$. 	 		
 		To introduce a control in \eqref{eq:particlemodel}, we replace the force $G$ depending on $T$ in \eqref{eq:particlemodel} by a force $F$ depending on a control $u\in U$. 
 		This results in an interacting particle model of the form
 		\begin{align}\label{eq:particlemodelu}
 		\frac{\di x_j}{\di t}=\frac{1}{N}\sum_{\substack{k=1\\k\neq j}}^N  F(x_j-x_k,u(x_j)).
 		\end{align}
 		A typical aspect of aggregation models is the competition of social interactions (repulsion and attraction) between particles and thus we assume that the  force $ F$  is of the form
 		\begin{align}\label{eq:totalforce}
 		F(d=d(x,y),u(x))= F_A(d,u(x))+ F_R(d,u(x)),
 		\end{align} 
 		where $d=d(x,y)=x-y$.  
 		Here,  $ F_R(d,u(x))$ denotes the repulsion force that a particle at location $y$ exerts on particle at location $x$ subject to the control parameter $u(x)$, and $ F_A$ is the attraction force a particle at location $y$ exerts on particle at location $x$, again subject to the control parameter $u(x)$. 
 		We assume that $F_R$ and $F_A$ are of the form
 		\begin{align}\label{eq:repulsionforce}
 		F_R(d=d(x,y),u(x))=\eta(x) f_R(\eta(x) |d|)d
 		\end{align}
 		and
 		\begin{align}\label{eq:attractionforce}
 		F_A(d=d(x,y),u(x))=\eta(x) f_A(\eta(x) |d|) R_{\theta(x)} \begin{pmatrix} 1 & 0 \\ 0 & \chi\end{pmatrix} R_{\theta(x)}^T d,
 		\end{align}
 		respectively, where the spatially inhomogeneous control parameter $u(x)$ is defined as $u(x) = (\theta(x), \eta(x))$. Here, we consider $\chi\in[0,1]$ and radially symmetric coefficient functions $f_R$ and $f_A$, where, again, $d=d(x,y)=x-y\in\R^2$. The rotation matrix $R_{\theta(x)}$ is defined as
 		\begin{align*}
 		R_{\theta(x)}=\begin{pmatrix} \cos(\theta(x)) &-\sin(\theta(x))\\
 		\sin(\theta(x)) & \cos(\theta(x))	\end{pmatrix}.
 		\end{align*}
 		The unusual form of the attraction force $F_A$ is motivated by \cite{patternformationanisotropicmodel,meanfieldfingerprints,stabilityanalysisanisotropicmodel,During2017,Merkel} where the direction of the interaction force depends on a spatially homogeneous or inhomogeneous tensor field $T=T(x)$ with
 		\begin{align*}
 		T(x):=\chi s(x)\otimes s(x) +l(x)\otimes l(x)\in\R^{2,2}
 		\end{align*}
 		for orthonormal  vector fields $s=s(x)$ and $l=l(x)\in\R^2$. Writing $s,l$ in polar coordinates results in 
 		$$s(x)=(-\sin(\theta(x)),\cos(\theta(x))), \quad l(x)=(\cos(\theta(x)),\sin(\theta(x)))$$
 		and the tensor field $T$ is given by 
 		$$T(x)=R_{\theta(x)} \begin{pmatrix} 1 & 0 \\ 0 & \chi\end{pmatrix} R_{\theta(x)}^T.$$	
 		This expression occurs in the definition of the attraction force $F_A$ in \eqref{eq:attractionforce}. 
 		
 		The parameter $\chi$ introduces an anisotropy to the force $F$ if $\chi<1$. The force $F$ along $l(x)=(\cos(\theta(x)),\sin(\theta(x)))$ is independent of $\chi$ and we assume that $F$ is short-range repulsive, long-range attractive along $l$. This implies that $F$ is also short-range repulsive, long-range attractive along $s(x)=(-\sin(\theta(x)),\cos(\theta(x)))$ for $\chi=1$, while for $\chi=0$ the total force $F$ along $s$ is purely repulsive. 
 		
 		For the forces $F_R$ and $F_A$, we make the following assumption:
 		\begin{ass}\label{ass:force}
 			We assume that the force coefficients $f_R$ and $f_A$ are bounded with 
 			\begin{align}\label{eq:forcecoeffperiodic}
 			f_R(|d|)\to 0,\quad f_A(|d|)\to 0\quad \text{ as } |d|\to 0.5.
 			\end{align}
 			Further, $f_R,f_A$ are continuously differentiable, implying
 			that the partial derivatives  of $F_R$ and $F_A$  are bounded, i.e.
 			\begin{align*}
 			\sup_{d\in\R^2}	\left|\frac{\partial F_R}{\partial d}\right|<+\infty, \quad \sup_{d\in\R^2} \left|	\frac{\partial F_A}{\partial d}\right|<+\infty
 			\end{align*}
 			and
 			\begin{align*}
 			\sup_{\eta\in [\eta_\text{min},\eta_\text{max}]}	\left|\frac{\partial F_R}{\partial \eta}\right|<+\infty, \quad \sup_{\eta \in[\eta_\text{min},\eta_\text{max}]} \left|	\frac{\partial F_A}{\partial \eta}\right|<+\infty.
 			\end{align*}
 		\end{ass}
 		Note that these conditions  are satisfied for the exponentially decaying force coefficients in \cite{patternformationanisotropicmodel,meanfieldfingerprints,stabilityanalysisanisotropicmodel,During2017,Merkel}. The rather unusual assumption \eqref{eq:forcecoeffperiodic} is considered to guarantee physically relevant forces through periodic extension on the torus which will be introduced in the following. 
 		
 		Motivated by the numerical simulations in \cite{patternformationanisotropicmodel,During2017}, the aim of this work is to study \eqref{eq:particlemodelu} on the torus and we consider $\Omega=\mathbb{T}^2\subset \R^2$. 
 		Since the torus $\T^2$ can be associated with the unit square $[0, 1]^2$ with periodic boundary conditions, it is useful to consider periodically defined forces for the associated discretized problems. Hence, we assume that $\bar F\colon\R^2\times U\to \R^2$ is the periodic extension of some force $F$, defined by
 		\begin{align}\label{eq:forceperiodic}
		\bar F(d+k,u(x))&= F(d,u(x)),\quad x\in\R^2,~ d=d(x,y)\in[-0.5,0.5]^2,~ k\in\Z^2,
 		\end{align}
 		for any $u\in U$,
 		see \cite{stabilityanalysisanisotropicmodel} for more details. 	As the  solutions are very sensitive to the scaling parameter $\eta,$ we restrict the space of controls $U$ to the space of admissible controls $\Uad \subset U$, defined as 
 		\begin{align}\label{eq:uad}
 		\Uad =\{u=(\theta,\eta) \in U \colon \eta \in [\eta_\text{min},\eta_\text{max}] \text{ a.s.}\}
 		\end{align}
 		 for  $0<\eta_\text{min}<\eta_\text{max}$. 
 		The  force $F$ and  its periodic extension $\bar F$ are Lipschitz continuous with respect to their  first variable with Lipschitz constant $C_d$, i.e.
 		$$|F(d,u(x))-F(\bar{d},u(x))|\leq C_d|d-\bar{d}|$$
 		for all $u\in\Uad$
 		where the Lipschitz constant $C_d$ is independent of $u=(\theta,\eta)$ due to the boundedness of $\eta$. Clearly, the partial derivatives with respect to $\theta$ are  bounded and hence, due to Assumption \ref{ass:force}, there exists a Lipschitz constant $C_u$ such that	
 		$$|F(d,u(x))-F(d,\bar{u}(x))|\leq C_u|u(x)-\bar{u}(x)|$$
 		for all $d\in\R^2$.

 		\subsection{State model}
 		Before formulating the state model on $\Omega =\mathbb{T}^2$, we introduce our notation. Let $\mathcal P(\domain)$ denote the space of Borel probability measures on $\domain$. By $\mathcal P_2(\Omega)$, we denote the space of Borel probability measures on $\domain$ with finite second moments, endowed with the 2-Wasserstein distance $W_2$, and by $\mathcal P_2^\emph{ac}(\Omega)\subset \mathcal P_2(\Omega)$ we denote the space of  Borel probability measures with finite second moments which are absolutely continuous with respect to the Lebesgue measure. For any map $g\colon \Omega\to \R^2$ and any measure $\mu\colon \Omega \to [0,+\infty]$ we denote by $g_\# \mu\colon \R^2 \to [0,+\infty]$ the pushforward measure, defined by $g_\# \mu(B)=\mu(g^{-1}(B))$ for any Borel set $B\subset \R^2$. In the following, we  consider the restriction of $U$ to $\Omega$, given by $U=\mathcal H(\Omega)\times \mathcal H(\Omega)$ and $\Uad$ as defined in \eqref{eq:uad}.
 		
 		Next, we formulate the state model. For this, we consider a control $u= (\theta,\eta)\in U$, i.e.\  $\theta(x) \in \R$ and $\eta(x) \in \R$ for all $x\in \R^2$, and introduce the interacting particle model 
 		\begin{align}\label{eq:particlemodelcontrol}
 		\frac{\di x_j}{\di t}=\frac{1}{N}\sum_{\substack{k=1\\k\neq j}}^N \bar F(x_j-x_k,u(x_j))
 		\end{align} 
 		for the periodic interaction force $\bar F$ on the unit torus $\Omega=\mathbb{T}^2\subset \R^2$.
 		The associated continuum model on  $\domain$ is given by 
 		\begin{align}\label{eq:macroscopiceq}
 		\begin{split}
 		\partial_t \rho(t,x)+\nabla_x\cdot \left[ \rho(t,x) (\bar F(\cdot,u(x)) \ast \rho(t,\cdot))(x)\right]&=0\qquad \text{in }[0,T]\times \domain
 		\end{split}
 		\end{align}  
 		for any $T>0$ and is equipped with initial data
 		\begin{align*}
 		\rho(0,\cdot) &=\rho^\text{in} \qquad \text{in } \domain
 		\end{align*}
 		for some given probability $\rho^\text{in}\in \mathcal P(\domain)$. 
 		
 		For $\bar F$ and $u$ given, \eqref{eq:macroscopiceq} has a unique global (weak measure) solution $\rho\in C([0, T ], \mathcal P_2(\Omega))$ for any $T>0$ which can be shown as in \cite[Thm 1.3.2]{Golse}. We refer to the weak solution of  \eqref{eq:macroscopiceq} at time $T>0$ as the solution to the state problem and we often write $\rho_t :=\rho(t,\cdot)$ for the solution of \eqref{eq:macroscopiceq} at time $t\in[0,T]$. 
  		We assume that $T>0$ is  sufficiently large so that the weak solution to \eqref{eq:macroscopiceq} at time $T$ is close to the corresponding steady state $\rho_\infty$ of \eqref{eq:macroscopiceq} satisfying
 		\begin{gather}\label{eq:stateproblem}
 		\rho_\infty(x) (\bar F(\cdot,u(x))\ast \rho_\infty)(x) = 0,  \qquad
 		\big(\bar F(\cdot,u(x))\ast \rho_\infty \big)(x) = \int_\Omega \bar F(x-y,u(x)) \rho_\infty(y) \di y.
  		\end{gather}
 	
 		Provided $\chi\geq 0$ is chosen sufficiently small,  patterns can be obtained as stationary  solutions to \eqref{eq:macroscopiceq} whose direction is controlled by the angle $\theta(x)$ for $x\in\Omega$. 
 		Since $s(x)$ rotates anticlockwise as $\theta(x)$ increases with $s(x)=(0,1)$ for $\theta(x)=0$,  $\theta(x)$ is the angle between the direction of the stationary  pattern at $x$ and the vertical axis. Note that rotations $\theta(x)+k\pi$ for any $k\in\Z$ result in the same direction of the pattern at $x$, implying that it is sufficient to consider $\theta(x)\in [0,\pi)$ for $x\in\Omega$. This can also be seen by the fact that we have $F_A(d,(\theta(x),\eta(x)))=F_A(d,(\theta(x)+\pi,\eta(x)))$ for any $\theta(x)\in[0,\pi)$ and hence also for its periodic extension. 	
 		
 		The solution of the mean-field PDE \eqref{eq:macroscopiceq} depends on the choice of initial data. In particular, this implies that the solution to \eqref{eq:stateproblem} is not unique. For concentrated initial data, a single vertical line along $s$ is expected as stationary solution, while for other initial data, 
 		more complex patterns may arise as stationary solution. The distance between those lines can be controlled by rescaling the total force which is controlled by the positive function $\eta$.

		\begin{rmrk}\label{rem:uniquenessstateproblem}
	   		The solution to  \eqref{eq:stateproblem} 	is not unique in general. To see this, note that  any solution $\rho_\infty$ of \eqref{eq:stateproblem} implies that $\rho_\infty +c$ for any constant $c \in \R$ is also a solution to \eqref{eq:stateproblem} since we have $\int_\Omega \bar F(x-y,u(x))\di y=0$ for any $x\in \Omega$. 
	   		To guarantee a unique solution, we consider an approximation of a stationary solution to \eqref{eq:stateproblem} as the state problem, given by the solution to the anisotropic aggregation equation \eqref{eq:macroscopiceq} after a sufficiently large time $T>0$  for specific initial data. 
	   	\end{rmrk}
	 	
	 	\begin{rmrk}\label{rem:propertyv}
	 			The term $(\bF(\cdot,u(x))\ast \rho)(x)$ for $\rho\in \mathcal{P}(\domain)$ in the state problem \eqref{eq:stateproblem} can be regarded as a macroscopic velocity field. We denote the space of Lipschitz continuous functions on $\domain$ by $\operatorname{Lip}(\domain)$ and we write $\langle \cdot,  \cdot \rangle$ for the scalar product on $\R^2$. The velocity field $v\colon \mathcal P(\domain)\times U_\emph{ad} \to \operatorname{Lip}(\domain), v(\rho,u)(x)=(\bF(\cdot,u(x))\ast \rho)(x)$ satisfies 
	 		\begin{subequations}\label{eq:vLipschitz}
	 		\begin{align}
	 		\langle v(\rho,u)(x)-v(\rho,u)(y),x-y \rangle&=\left\langle \int_\Omega (\bF(x-z,u(x))-\bF(y-z,u(y)))\rho(z)\di z,x-y\right\rangle \\
	 		& \leq C_d |x-y|^2, \quad x,y\in\Omega
	 		\end{align}
	 		\end{subequations}
	 		for all $(\rho,u)\in \mathcal P(\domain)\times \Uad$
	 		where the Lipschitz constant $C_{d}>0$  is independent of $(\rho,u).$ 
	 		
	 		For any function $g$ on $\domain$, we write $\|g\|_\infty=\sup_{x\in \domain}|g(x)|$ for the supremum norm and we have for any  $(\rho,u),(\bar{\rho},\bar{u})\in \mathcal P(\domain)\times U_\emph{ad}$:
	 		\begin{align*}
	 		\|v(\rho,u)-v(\bar{\rho},\bar{u})\|_\infty\leq \|v(\rho,u)-v(\bar{\rho},u)\|_\infty+\|v(\bar{\rho},u)-v(\bar{\rho},\bar{u})\|_\infty.
	 		\end{align*}
	 		Indeed, we have 
	 		$$\| v(\rho,u)-v(\bar{\rho},u)\|_\infty\leq C_{v} W_2(\rho,\bar{\rho})$$
	 		for some constant $C_{v}>0$ 
	 		and
	 		$$\|v(\bar{\rho},u)-v(\bar{\rho},\bar{u})\|_\infty=\sup_{x\in\Omega}\left|\int_\Omega (\bF(x-z,u(x))-\bF(x-z,\bar{u}(x)))\bar{\rho}(z)\di z\right| \leq C_u \|u-\bar{u}\|_\infty$$
	 		due to the continuous embedding of $\mathcal H(\Omega)$ in $L^\infty(\Omega)$.
	 		This implies
	 		\begin{align}\label{eq:vsup}
	 		\|v(\rho,u)-v(\bar{\rho},\bar{u})\|_\infty\leq C_{v} W_2(\rho,\bar{\rho})+C_u \|u-\bar{u}\|_\infty.
	 		\end{align}
	 		For weak solutions $\rho, \bar{\rho}$ of \eqref{eq:macroscopiceq} with initial data $\rho^\text{in}, \bar{\rho}^\text{in}$ and controls $u,\bar{u}$, respectively, one can show  that there exists positive constants $a,b$ such that
	 		\begin{align}\label{eq:W2bound}
	 		W_2^2(\rho_t,\bar{\rho}_t)\leq \left(W_2^2(\rho^\text{in},\bar{\rho}^\text{in}) +b \|u-\bar{u}\|_\infty^2\right)\exp(at)
	 		\end{align}
	 		for all $t\in[0,T]$ where we used again the continuous embedding of $\mathcal H(\Omega)$ in $L^\infty(\Omega)$. The proof is based on computing the derivative of $W_2^2(\rho_t,\bar{\rho}_t)$ with respect to $t$, estimates using  the inequalities \eqref{eq:vLipschitz}, \eqref{eq:vsup}, and the Gronwall inequality, see \cite{BPTT} for more details.
	 	\end{rmrk}

	 	\subsection{Optimal control problem in the macroscopic setting}\label{sec:optimalcontrolmacro}
	 	With the state problem defined in \eqref{eq:stateproblem}, we can now formulate the associated optimal control problem. For $\bF$ and $u$ given, we identify the solution of the state problem as the weak solution of the macroscopic problem  \eqref{eq:macroscopiceq} at time $T>0$ and denote it by $\rho$ in the following.  	
		Given the distribution function $\rhod \in \mathcal{P}(\domain)$ of a desired fingerprint pattern, the task at hand is to find the force field $\bF$ characterized by $u = (\theta, \eta)$ corresponding to the given density $\rhod$. This task can be regarded as an inverse problem which is mathematically formulated as an optimal control problem with a PDE constraint. 
	We  define the cost functional
	\begin{align}\label{eq:costfunctional}
	 \J(\rho_T,u) =  \frac{1}{2} \mathcal W_2^2(\rho_T,\rho_\text{des})   + \frac{\lambda_1}{2} \|\theta - \theta_\text{ref}\|^2_{\mathcal H(\Omega)} + \frac{\lambda_2}{2} \|\eta - \eta_\text{ref}\|^2_{\mathcal H(\Omega)}
	\end{align}
	where $\lambda_1, \lambda_2 > 0$ are parameters, and $\theta_\text{ref}, \eta_\text{ref}\in \mathcal H(\Omega)$ are given reference values. To summarize, we obtain:
	\begin{prblm}\label{prob:macroscopic}
		Find $u\in U_\emph{ad}$ such that
	\begin{align*}
		(\rho_T(u),u)=\argmin_{\rho_T,u} \J(\rho_T,u) \quad \text{subject to }\eqref{eq:macroscopiceq}.
	\end{align*}
	\end{prblm}
	Considering the well-defined solution operator $S\colon U \to \mathcal P(\domain)$ with $Su= \rho_T$ associated with the unique weak solution of  \eqref{eq:macroscopiceq}, the optimization problem is then given by
	$$\min \tilde{\J}(u)$$ where $\tilde{\J}$ is defined by $\tilde \J(u) :=\J(Su,u)$ and is called the reduced cost functional.

	\section{First-order optimality conditions in the macroscopic setting}\label{sec:foc}
	The main objective of this section is to derive the first-order optimality conditions (FOC) for the optimal control problem by using the Lagrangian approach based on Wasserstein calculus. The arguments of this section are similar to the ones in \cite{BPTT}. Since  there exists a unique global (weak measure) solution $\rho\in C([0, T ], \mathcal P_2(\Omega))$ of the mean-field PDE \eqref{eq:macroscopiceq} for any $u\in U$, 
	we define the state operator 
	\begin{equation}\label{eq:weak}
	E(\rho,u) [\varphi] := \langle \varphi_T,  \rho_T \rangle - \langle \varphi_0, \rho^\text{in} \rangle - \int_0^T \langle \partial_t \varphi(t,x) +(\bF(\cdot,u(x) ) \ast \rho_t)(x) \cdot \nabla \varphi(t,x) , \rho(t,x) \rangle \di t 
	\end{equation}
	with $E(\rho,u)=0$
	for all $\varphi\in \mathcal A:= C_c^1([0,T]\times\Omega)$. 
	Further, let
	\begin{align}\label{eq:jfunctionals}
	 \J_1(\rho) :=  \frac{1}{2} \mathcal W_2^2(\rho_T,\rho_\text{des}),\quad \J_2(u) :=   \frac{\lambda_1}{2} \|\theta - \theta_\text{ref}\|^2_{\mathcal H(\Omega)},\quad \J_3(u) :=  \frac{\lambda_2}{2} \|\eta - \eta_\text{ref}\|_{\mathcal H(\Omega)}^2. 
	 \end{align}
	Hence, we have $\J(\rho, u) = \J_1(\rho) + \J_2(u)+\J_3(u)$. Since the functional $\J$ is not handy for deriving the first-order optimality conditions, we consider the extended Lagrangian $\mathcal I$ defined by
	$$ \min_{(\rho, u)} \mathcal I(\rho,u) = \min_u \Big\{ \J_2(u)+\J_3(u) + \min_{\rho} \sup_{\varphi \in \mathcal A} \big\{ \J_1(\rho) + E(\rho,u)[\varphi] \big\} \Big\} = \min_u \Big\{ \J_2(u)+\J_3(u) + \kappa(u) \Big\}$$
	with $$ \kappa(u) = \min_{\rho} \sup_{\varphi \in \mathcal A} \Big\{ \J_1(\rho) + E(\rho,u)[\varphi] \Big\}. $$
	Note that $E$ has the property $\sup_{\varphi \in \mathcal A} E(\rho,u)[\varphi] \ge 0$ as $\varphi \equiv 0$ implies $E(\rho,u)[0] = 0$ for every $(\rho,u).$ Therefore, if $E(\rho,u)[\varphi] >0$ for some $\varphi$, the linearity in $\varphi$ of $E$ yields $E(\rho,u)[a\varphi] = a E(\rho,u)[\varphi]$ for every $a > 0$ which shows that $\sup_\varphi E(\rho,u)[\varphi] = + \infty.$ 
	
	We say a pair $(\rho,u)\in C([0, T ], \mathcal P_2(\Omega))\times U$ is admissible if $E(\rho, u)[\varphi] = 0$ for all $\varphi \in \mathcal A$. 
	In the following, we derive a necessary condition for an admissible pair $(\rho,u)$ to be a stationary point. For this, we consider perturbations of the admissible pair $(\rho,  u)$, given by an admissible perturbation $u^\delta=u+ \delta h \in U$   of $ u$ for $\delta\geq 0$ and   $h=(\theta,\eta) \in U$, 
	and the associated  $\rho^\delta \in \mathcal P_2(\Omega)$ satisfying $E(\rho^\delta, u^\delta)[\varphi] = 0$ for all $\varphi \in \mathcal A$. 	
	  These conditions result in the following assumption:
	\begin{ass}\label{ass:stationarypoint}
		Let $(\rho,  u)$ be an admissible  pair and suppose that $h\in U$.
		Suppose $\bar{\delta}>0$ is given such that for any $0\leq \delta \leq \bar{\delta}$ the perturbation $u^\delta :=  u + \delta h$ satisfies
			\begin{itemize}
			\item $u^\delta $ is an admissible control, i.e.\ $u^\delta \in U$.
			\item There exists   $\rho^\delta \in  C([0, T ], \mathcal P_2(\Omega))$ such that $E(\rho^\delta, u^\delta)[\varphi] = 0$ for all $\varphi \in \mathcal A$. 
		\end{itemize}
	\end{ass}		
	For an admissible pair $(\rho,u)$ satisfying Assumption \ref{ass:stationarypoint} we have
	\begin{align*}
	\kappa(u^\delta) &= \min_{\rho}\sup_{\varphi\in\mathcal A}\, \Bigl\{  \J_1(\rho) + E(\rho,u^\delta)[\varphi]\Bigr\} = \J_1(\rho^\delta) \\
	&= \J_1(\rho^\delta) - \J_1(\rho) + \min_{\rho}\sup_{\varphi\in\mathcal A}\, \Bigl\{  \J_1(\rho) + E(\rho, u)[\varphi]\Bigr\} \\
	&= \J_1(\rho^\delta) - \J_1(\rho) + \kappa( u),
	\end{align*}
	and the directional derivative of $\mathcal G := \J_2 + \J_3+\kappa$ at $ u$ along $h$ is given by
	\begin{align*}
	\lim_{\delta\to 0}\frac{\mathcal G(u^\delta)-\mathcal G( u)}{\delta} = \lim_{\delta\to 0}\frac{[\J_1(\rho^\delta) - \J_1(\rho)] + [\J_2(u^\delta) - \J_2( u)]+ [\J_3(u^\delta) - \J_3( u)]}{\delta},
	\end{align*}
	which depends on the relationship between $\rho^\delta$ and $\rho$. 
	
	\begin{rmrk}
	Using estimate \eqref{eq:W2bound}, we obtain
		\begin{align*}
			W_2(\rho^\delta_t,\rho_t)\leq   \delta \sqrt{b} \|h\|_\infty\exp\left(\frac{aT}{2}\right)
		\end{align*}	
		for any $t\in[0,T]$.
		Hence, for $t\in[0,T]$ given, the curve $[0,\infty)\ni \delta \mapsto \rho^\delta_t \in \mathcal P_2(\Omega)$ which starts from $\rho^0_t$ at $\delta=0$ is absolutely continuous with respect to the 2-Wasserstein distance and there exists $\psi_t\in L^2(\rho_t;\Omega)$  satisfying \cite[Proposition~8.4.6]{Ambrosio}
		\begin{equation}\label{eq:vector-field}
		\lim_{\delta\to 0} \frac{W_2(\rho^\delta_t,(id+\delta \psi_t)_\#\rho_t)}{\delta} = 0.
		\end{equation}
		Furthermore,
		\begin{align*}
		W_2^2((id + \delta\psi_t)_\#\rho_t,\rho_t) \le \iint_{\Omega\times\Omega} |x + \delta\psi_t(x)-x|^2\di \rho_t(x) = \delta^2\iint_{\Omega\times\Omega} |\psi_t(x)|^2\di \rho_t(x).
		\end{align*}
		In particular, we have that
		\begin{align*}
		\limsup_{\delta\to 0} \frac{W_2(\rho^\delta_t,\rho_t)}{\delta} = \limsup_{\delta\to 0} \frac{W_2((id+\delta \psi_t)_\#\rho_t,\rho_t)}{\delta}\le \sqrt{\iint_{\Omega\times\Omega} |\psi_t|^2\dd\rho_t(x)}.
		\end{align*}
	\end{rmrk}

In the following we establish an explicit relationship between the perturbations $\psi_t$ and $h$ as in \cite[Thm 3.4]{BPTT}. Note that we omit the proofs of the following lemma and the theorem because they are very similar to the ones in the reference.

\begin{lmm}\label{lem:4.4}
	Let $(\rho,\con)$ be an admissible pair and let $\delta>0$ sufficiently small, $\h\in U$ and $\con^\delta=\con + \delta {\h}$ such that 
	\begin{enumerate}
		\item[(i)] $\con^\delta\in U$, and
		\item[(ii)] there exists $\rho^\delta\in \mathcal C([0,T],\mathcal P_2^\emph{ac}(\Omega))$ satisfying $E(\rho^\delta,\con^\delta)=0$.
	\end{enumerate}
	Suppose that $\psi\in \mathcal C^1((0,T)\times \Omega)$ with $\psi_0\equiv 0$ satisfies
	\begin{align}\label{eq:linearization}
	\partial_t\psi_t +  D\psi_t\,v(\rho_t,\con) = \calK(\rho_t,\con)[\psi_t,\h]\qquad \text{for $\rho_t$-almost every $x\in\Omega$},
	\end{align}
	where  $(t,x)\mapsto \calK(\rho_t,\con)[\psi_t,\h](x)$ is  a bounded Borel map
	satisfying
	\begin{align}\label{eq:K_limit}
	\lim_{\delta\to 0} \int_0^T \int_\Omega \left|\frac{v(\nu_t^\delta,\con^\delta)\circ(id + \delta\psi_t)(x) - v(\rho_t,\con)(x)}{\delta} -  \calK(\rho_t,\con)[\psi_t,\h](x)\right|^2 \dd\rho_t(x)\dd t = 0
	\end{align}
	with $\nu_t^\delta := (id + \delta\psi_t)_\# \rho_t$.
 Then, \eqref{eq:vector-field} holds for this $\psi$, i.e.\
	\[
	\lim_{\delta\to 0} \frac{W_2(\mu_t^\delta,(id+\delta \psi_t)_\#\mu_t)}{\delta} = 0.
	\]
\end{lmm}

\begin{rmrk}
	Note that for any $h\in U$ we obtain by Taylor expansion:
	\begin{equation}\label{eq:calK}
	\calK(\rho,\con)[\psi,\h] = Dv(\rho,\con)\psi  - \int_\Omega (\nabla_{x} \bF)(\cdot -y,u)\psi(y)\dd  \rho(y)-\nabla_u v(\rho,u) h+\mathcal O(\delta),
	\end{equation}
	and \eqref{eq:linearization} may be written as
	\[
	\partial_t\psi_t +  \{\psi_t,v(\rho_t,\con)\} =  -\int_\Omega (\nabla_{x} \bF)(\cdot -y,u)\psi_t(y)\dd  \rho_t(y)-\nabla_u v(\rho_t,u) h ,
	\]
	where $\{\cdot,\cdot\}$ denotes the Lie bracket given by $\{\xi,\psi\}=(D\xi) \psi - (D\psi) \xi$ for vector fields $\xi,\psi$.
\end{rmrk}

The existence of $\psi \in C_b^1((0,T) \times \Omega)$ satisfying the assumptions of Lemma~\ref{lem:4.4} is provided in the following statement.
\begin{thrm}
Suppose that the assumptions of Lemma~\ref{lem:4.4} hold. For the velocity field $v \colon \mathcal P_2(\Omega) \times U \rightarrow \emph{Lip}_\emph{loc}(\Omega)$ given by $v(\rho,\con)(x) =  (F(\cdot,\con(x)) \ast \rho)(x)$ there exists $\psi \in \mathcal C_b^1((0,T) \times \Omega)$ with $\psi_0 = 0$ satisfying
\[
\partial_t \psi_t + D\psi_t v(\rho_t,\con) = \calK(\rho_t,\con)[\Psi_t,h] \quad \text{for $\rho_t dt$-almost every\ } (t,x) \in (0,T) \times \Omega, 
\]
where $\calK$ is given by \eqref{eq:calK}.
\end{thrm}

Now, we are able to state the first-order necessary condition for $(\rho,\con)$ to be a stationary point. For this, the Gâteaux-derivatives of $\J_1,\J_2$ and $\J_3$ in \eqref{eq:jfunctionals} are required. By \cite[Proposition~8.5.2]{Ambrosio}, we have
	\begin{align}\label{eq:deltaJ_1}
\delta_{\rho} \J_1(\rho_T)(x) = t_{\rho_T}^{\rho_\text{des}} (x) -x
\end{align}
	where $ t_{\rho_T}^{\rho_\text{des}} (x)$ is the unique optimal transport map from $\rho_T$ to $\rho_\text{des}$. 
For $u=( \theta, \eta)\in U$ and $h=(h_\theta,h_\eta)\in U$, we have
\begin{align}\label{eq:deltaJ_2}
\dd \J_2( u)[h] +	\dd \J_3(\con)[\h]= \lambda_1  \langle h_\theta,\theta- \theta_\text{ref}\rangle_{\mathcal H(\Omega)}  + \lambda_2  \langle h_\eta,\eta- \eta_\text{ref}\rangle_{\mathcal H(\Omega)}
\end{align}
for the inner product $\langle \cdot,\cdot \rangle_{\mathcal H(\Omega)}$ on $\mathcal H(\Omega)$
since the  Gâteaux-derivatives are given by
\begin{align*}
\dd \|\theta- \theta_\text{ref}\|_{H(\Omega)}[h_\theta]=\left\langle h_\theta,\frac{\theta- \theta_\text{ref}}{\|\theta - \theta_\text{ref}\|_{\mathcal H(\Omega)}} \right\rangle_{\mathcal H(\Omega)} \quad \text{for} \quad \theta- \theta_\text{ref}\neq 0,\\
\dd \|\eta - \eta_\text{ref}\|_{\mathcal H(\Omega)}[h_\eta]=\left\langle h_\eta,\frac{\eta- \eta_\text{ref}}{\|\eta - \eta_\text{ref}\|_{\mathcal H(\Omega)}} \right\rangle_{\mathcal H(\Omega)}\quad \text{for} \quad \eta- \eta_\text{ref}\neq 0.
\end{align*}
\begin{thrm}\label{thm:optW}
	Let $(\bar\rho,\bar\con)$ be an optimal pair, $\J_1,\J_2,\J_3$ be G\^ateaux-differentiable. Suppose that $\h \in U$ and assume that there exists $\bar \delta>0$ such that for any $0\leq \delta \leq \bar \delta$ it holds $\bar u +\delta h\in U$. Then,
		\begin{align}\label{eq:optimality}
	\dd \J_2(\bar\con)[\h]+	\dd \J_3(\bar\con)[\h] + \int_\Omega \langle \delta_\rho \J_1(\bar\rho_T), \psi_T\rangle\dd \bar \rho_T= \lim\limits_{\delta \rightarrow 0} \frac{\mathcal G(\bar u+\delta h) - \mathcal G(\bar u)}{\delta} = 0
	\end{align}
	where $\J_1$, $\J_2$, $\J_3$  are defined in \eqref{eq:deltaJ_1}, \eqref{eq:deltaJ_2} and $t\mapsto\psi_t\in L^2(\bar\rho_t)$ satisfies \eqref{eq:linearization} with $\psi_0=0$. 
\end{thrm}

The optimality condition \eqref{eq:optimality} together with \eqref{eq:deltaJ_1} and \eqref{eq:deltaJ_2} is implicit. In order to formulate an optimization algorithm, we derive an explicit form of the optimality conditions by computing the adjoint.
Thus, to derive the adjoint-based first-order optimality system, we begin with the dual problem corresponding to \eqref{eq:linearization}. It can be obtained by testing  \eqref{eq:linearization} with a family of vector-valued measures $(m_t)_{t\in(0,T)}$ with $m_t\in \mathcal{P}(\Omega;\R^2)$, resulting in 
\begin{align*}
\int_0^T \!\!\int_\Omega \Bigl(\partial_t\psi_t + D\psi_t\,v(\bar\rho_t,\bar\con)  - \calK(\bar\rho_t,\bar\con)[\psi_t,h]\Bigr)\cdot \dd m_t \dd t = 0
\end{align*}
or equivalently 
\begin{align*}
\int_0^T\int_\Omega \left(\partial_t\psi_t +  \{\psi_t,v(\bar \rho_t,\bar \con)\}   +\int_\Omega (\nabla_{x} \bF)(\cdot -y,\bar u)\psi_t(y)\dd \bar \rho_t(y) + \nabla_u v(\bar\rho_t,\bar\con) h \right)\cdot \di m_t \dd t=0.
\end{align*}
Integrating by parts and using $\psi_0=0$, we obtain 
\begin{align*}
&\int_0^T \langle  \partial_t m_t  + \nabla \cdot (v(\bar \rho_t, \bar u) \otimes  m_t) + \nabla v(\bar \rho_t ,\bar u)  m_t -\bar \rho_t\int_\Omega (\nabla_{x} \bF)(y-\cdot,\bar u)\dd  m_t(y),  \psi_t \rangle \dd t\\&=\int_\Omega \psi_T \cdot \di m_T+\int_0^T \int_\Omega \nabla_u v(\bar \rho_t,\bar u) h\cdot \di m_t \di t.
\end{align*}
We choose $\bar m_t$ such that  the dual problem
\[\partial_t \bar m_t + \nabla \cdot (v(\bar \rho_t, \bar u) \otimes \bar m_t) + \nabla v(\bar \rho_t ,\bar u) \bar m_t - \bar \rho_t \int_\Omega (\nabla_x \bF)(y-\cdot,\bar u) \dd \bar m_t(y) = 0
\]
subject to the terminal condition $\bar m_T = \bar \rho_T \delta_\rho \J_1(\bar \rho_T)$. By using the optimality condition \eqref{eq:optimality}, we obtain
\begin{align*}
-\dd \J_2(\bar\con)[\h]-\dd \J_3(\bar\con)[\h] +\int_0^T \int_\Omega \nabla_u v(\bar \rho_t,\bar u) h\cdot \di \bar m_t \di t=0
\end{align*}
for any $h\in U$. The optimality condition can be summarised as:
\begin{thrm}
A minimizing pair $(\bar \rho,\bar u)$ satisfies the state problem
\begin{align*}
\partial_t \bar\rho_t + \nabla \cdot (\bar \rho_t v(\bar\rho_t,\bar u)) &= 0
\end{align*} 
subject to the initial condition $\bar \rho(0,\cdot) = \bar \rho_0=\rho^\text{in},$
and the optimality condition
\begin{align}\label{eq:optimality2}
\dd \J_2(\bar u) +\dd \J_3(\bar u)= \int_0^T\int_\Omega \nabla_u v(\bar \rho_t, \bar u)  \dd \bar m_t \dd t,
\end{align}
where $\bar m_t$ solves the adjoint equation given by
\[
\partial_t \bar m_t + \nabla \cdot (v(\bar \rho_t, \bar u) \otimes \bar m_t) + \nabla v(\bar \rho_t ,\bar u) \bar m_t - \bar \rho_t \int_\Omega (\nabla_x \bF)(y-\cdot,\bar u) \dd \bar m_t(y) = 0
\]
subject to the terminal condition $\bar m_T = \bar \rho_T \delta_\rho \J_1(\bar \rho_T).$
\end{thrm}

\section{Discretization of the optimality conditions}\label{sec:focdiscrete}

In this section, we formally derive the adjoints and the optimality conditions for the discrete optimal control problem. In addition, we introduce the discretised reduced cost functional and its gradient which are required for applying adjoint-based descent methods to solve the discretised control problems numerically. 
\subsection{Discrete adjoint system}
We  derive the discrete system in order to formulate an algorithm for numerical simulations. Note that the domain  $\Omega=\mathbb{T}^2$ can be regarded as the
unit square $[0, 1]^2$ with periodic boundary conditions. First, we make an ansatz, similar to \cite{BPTT}, to obtain an equation for the vector-valued adjoint variable. Indeed, assuming $|\bar m_t| \ll \bar \rho_t$ for every $t \in [0,T]$ yields the existence of a vector field $\bar \xi_t \colon \Omega\to \R^2$ such that $\bar m_t = \bar \xi_t \bar \rho_t$ where $\bar \rho_t$ is the weak solution to the state equation \eqref{eq:macroscopiceq}. The dual problem allows us to obtain an equation for $\bar \xi_t$ given by
\[
\partial_t \bar \xi_t + D \bar\xi_t  v(\bar \rho_t, \bar u)= - \nabla v(\bar \rho_t ,\bar u) \bar \xi_t + \int_\Omega (\nabla_x \bF)(y-\cdot,\bar u)  \xi_t(y)\dd \bar \rho_t(y), 
\]
implying that we can write the adjoint equation as
\begin{align}\label{eq:adjointcont}
\partial_t \bar\xi_t + \nabla(v(\bar \rho_t, \bar u) \cdot \bar \xi_t) =  \int_\Omega  (\nabla_x \bF)(y-\cdot,\bar u) \bar \xi_t(y) \dd \bar \rho_t(y).
\end{align}
This PDE is equipped with the terminal condition $\bar \xi_T = \delta_\rho \J_1(\bar \rho_T) $ where $\delta_\rho \J_1(\bar \rho_T)(x) =  t_{\rho_T}^{\rho_\text{des}} (x) -x $ by \eqref{eq:deltaJ_1}. Note that \eqref{eq:adjointcont} has a unique solution $\xi \in C([0,T]\times \Omega)$ such that $\xi(t)\in \text{Lip}_b(\Omega)$ for $t\in[0,T]$. We omit the proof, as it is very similar to \cite[Thm 3.10]{BPTT}. Using $\bar m_t = \bar \xi_t \bar \rho_t$,  \eqref{eq:optimality2} can be written as
\begin{align}\label{eq:optimality3}
\begin{split}
\dd \J_2(\bar u) +\dd \J_3(\bar u)&= \int_0^T\int_\Omega \nabla_u v(\bar \rho_t, \bar u) \bar \xi_t \dd \bar \rho_t \dd t\\&=\int_0^T\int_\Omega \int_\Omega \nabla_u \bar F(x-y, \bar u(x)) \bar \xi_t(x) \dd \bar \rho_t(y)\dd \bar \rho_t(x) \dd t.
\end{split}
\end{align}

Now, we formally derive the associated optimality conditions on the particle level. For this, we consider the equation of characteristics for the mean-field PDE  \eqref{eq:macroscopiceq}, given by 
\[
\frac{\dd x }{\dd t}  = v( \rho^N, \bar u)(x)=\int_\Omega \bF(x - y,\bar u(x)) \dd\rho^N(y)
\]
 with the empirical measure 
\[
\rho^N(t,x) = \frac{1}{N}\sum_{j=1}^N \delta_0(x-x_j(t))
\]
of particle positions  $x_i\in\R^2,\; i=1,\dots,N$. The periodicity of the boundary conditions in the macroscopic setting is mimicked in the discrete setting by considering periodic forces defined by \eqref{eq:forceperiodic}. We emphasize that this is crucial, as otherwise the same two particles may interact various times. For a given control $u\in U$, this leads to the particle system 
\begin{align}\label{eq:discretemodel}
\frac{\dd x_i }{\dd t} = \frac{1}{N} \sum_{j=1}^N \bF(x_i - x_j, u(x_i)), \quad i=1,\dots,N,
\end{align}
subject to the initial conditions
\begin{align}\label{eq:discreteinitial}
 x_i (0)=X_i, \quad i=1,\dots,N,
\end{align}
 where  $X_i\in \Omega,\; i=1,\dots,N,$ are given independent realizations of random variables  with $\text{law}(X_i) = \rho_0.$ Note that \eqref{eq:discretemodel} is identical to \eqref{eq:particlemodelcontrol}.
 
 We introduce adjoint variables $\xi_i=\xi(x_i)\in\R^2,\; i=1,\dots,N$, and
 the discrete adjoint system reads
\begin{align}\label{eq:adjointdiscrete}
\frac{\dd}{\dd t} \xi_i = \frac{1}{N}\sum_{j=1}^N \nabla_x \bF(x_i - x_j, u(x_i)) \xi_i - \frac{1}{N} \sum_{j=1}^N \nabla_x \bF(x_i-x_j, u(x_i)) \xi_j 
\end{align}
with terminal condition $\xi_i(T) =  \delta_\rho \J_1(x_i(T))=t_{\rho_T}^{\rho_\text{des}} (x_i(T)) -x_i(T).$ 

\subsection{First-order optimality conditions}
We introduce the cost functional in the discrete setting. For this, we replace the probability measure $\rho_\text{des}$  by the empirical measure $$\rho_\text{des}^N=\frac{1}{N}\sum_{i=1}^N \delta_{x_i^\text{des}},$$
where $x_i^\text{des}\in \Omega, \; i=1,\dots,N$, are assumed to be given.
Then, the discrete analogue $\J^N\colon \R^{2N}\times U \to \R$ to the cost functional \eqref{eq:costfunctional} in the continuum setting is given by
\begin{align}\label{eq:costfunctionaldiscrete}
\begin{split}
&\J^N((x_1(T),\ldots,x_N(T)),u) \\&= \frac{N}{2} \mathcal W_2^2\left(\frac{1}{N}\sum_{i=1}^N \delta_{x_i(T)},\rho_\text{des}^N\right)   + \frac{\lambda_1}{2} \|\theta - \theta_\text{ref}\|^2_\infty+ \frac{\lambda_2}{2} \|\eta - \eta_\text{ref}\|^2_\infty
\end{split}
\end{align}
where $(\theta_\text{ref},\eta_\text{ref})\in\Uad$ is assumed to be given reference data. 
\begin{rmrk}\label{rem:rescaling}
	Note that we rescale the first term of the cost functional by $N$. This is necessary to have well-balanced terms in the discrete Lagrangian below, see \eqref{eq:discreteLagrangian}. Indeed, as $N \to \infty$ we obtain infinity many terms in the dual  accounting for the constraint, whereas the Wasserstein distance has fixed order. This scaling is also reported in \cite{FornasierScaling,BPTT}.
\end{rmrk}
Denoting $x^N=(x_1,\ldots,x_N)$, the discrete optimal control problem is given by
\begin{prblm}\label{prob:discrete}
	For $N\in \mathbb{N}$ fixed, find $u^N\in U_\emph{ad}$ such that
	\begin{align*}
	(x^N(T),u^N)=\argmin_{x^N(T),u^N} \mathcal{J}^N (x^N(T),u^N) \quad \text{subject to }\eqref{eq:discretemodel},
	\end{align*}
	where $x^N=x^N(u^N)$.
\end{prblm}
	We define the state operator  $e^N$ by 
	\begin{align*}
	e^N(x,u)=\begin{pmatrix}
	\frac{\dd}{\dd t} x(t) - \frac{1}{N} \left(\sum_{j=1}^N \bF(x_i(t) - x_j(t), u(x_i(t)))\right)_{i=1}^N  \\
	x(0)-X
	\end{pmatrix}
	\end{align*}
	and the discrete state system \eqref{eq:discretemodel}--\eqref{eq:discreteinitial} can be rewritten as $e^N(x,u)=0$. 
	Its weak formulation is given by
	\begin{align*}
	&\langle e^N(x,u),(\xi,\zeta)\rangle\\&= \int_0^T \left(
	\frac{\dd}{\dd t} x(t) - \frac{1}{N} \left(\sum_{j=1}^N \bF(x_i(t) - x_j(t), u(x_i(t)))\right)_{i=1}^N \right) \cdot\xi(t)\di t+
	(x(0)-X) \cdot \zeta\\&=\sum_{i=1}^N\int_0^T \left(
	\frac{\dd}{\dd t} x_i(t) - \frac{1}{N} \left(\sum_{j=1}^N \bF\left(x_i(t) - x_j(t), u(x_i(t))\right)\right) \right)\cdot \xi_i(t)\di t+
	(x(0)-X) \cdot \zeta
	\end{align*}
	where $\langle \cdot,\cdot \rangle$ denotes the scalar product on $\R^2$.
	For Lagrange multipliers $(\xi,\zeta)$, the Lagrangian corresponding to Problem \ref{prob:discrete} reads
	\begin{equation}\label{eq:discreteLagrangian}
	\mathcal L(x,u,\eta,\zeta)=\J^N((x_1(T),\ldots,x_N(T)),u)+\langle e^N(x,u),(\xi,\zeta)\rangle
	\end{equation}
	for $N\in\mathbb{N}$ fixed. 
	We have 
	\begin{align*}
	\dd_u\J^N (x, u)[h] = \dd \J_2(u)[h] +\dd \J_3(u)[h] 
	\end{align*}
	and
	\begin{align*}
	\langle \di_u e^N(x,u)[h],(\xi,\zeta)\rangle = - \frac{1}{N} \sum_{i=1}^N \sum_{j=1}^N \int_0^T 
	 \nabla_u \bar F(x_i(t) - x_j(t), u(x_i(t))) h(x_i(t))\cdot \xi_i(t)\di t,
	\end{align*}
	implying
	\begin{align*}
	 \dd \J_2(u)[h] +\dd \J_3(u)[h] 
	 - \frac{1}{N^2} \sum_{i=1}^N \sum_{j=1}^N \int_0^T 
	\nabla_u \bar F(x_i(t) - x_j(t), u(x_i(t))) h(x_i(t))\cdot \xi_i(t)\di t=0
	\end{align*}
	for any $h\in U$.
	The first-order optimality condition in the rescaled (see Remark~\ref{rem:rescaling}) microscopic setting reads:
\begin{thrm}
	Let $N\in \mathbb N$ be given and let $(x^N,u^N)$ be an optimal pair. The optimality condition corresponding to Problem \ref{prob:discrete} reads
	\begin{align}\label{eq:optimalitydiscrete}
	\dd \J_2(u^N)+\dd \J_3(u^N)
	=\frac{1}{N^2} \sum_{i=1}^N \sum_{j=1}^N \int_0^T 
	\nabla_u \bF(x_i(t) - x_j(t), u^N(x_i(t)))  \xi_i(t)\di t
	\end{align}
	where $\xi$ satisfies \eqref{eq:adjointdiscrete} with terminal condition $$\xi_i(T) =  N\delta_\rho \J_1(x_i(T))=N (t_{\rho_T}^{\rho_\text{des}} (x_i(T)) -x_i(T)).$$
\end{thrm}

\subsection{Gradient of the reduced cost functional}
In this section, we introduce the discretised reduced cost functional and its gradient. Motivated by the control-to-state operator $S$ and the reduced cost functional $\tilde \J$ in the mean-field setting, introduced in Section \ref{sec:optimalcontrolmacro}, we consider the control-to-state operator $S^N\colon U \to \R^{2N};\  u \mapsto (x_1,\ldots,x_N)$ where $(x_1,\ldots,x_N)$ satisfies the forward particle system \eqref{eq:discretemodel} on $[0,T]$ with initial conditions \eqref{eq:discreteinitial}.
The reduced cost functional $\tilde \J^N$ in the discrete setting is then given by $\tilde \J^N(u):= \J^N(S^N(u),u)$ where $\J^N$ is defined in \eqref{eq:costfunctionaldiscrete}. 
Since the force $\bF$ satisfies  Assumption~\ref{ass:force}, we can implicitly obtain $\di S^N (u)$ via 
\begin{align*}
0=\di_x e^N(S^N (u),u) [\di S^N (u)]+ \di_u e^N(S^N (u),u),
\end{align*}
i.e.\
\begin{align*}
\di S^N (u)[h]=-(\di_x e^N(S^N (u),u))^{-1} \di_u e^N(S^N (u),u)[h].
\end{align*}
The Gâteaux derivative $\tilde \J^N$ in the direction $h=(h_\theta, h_\eta)\in U$ is obtained from
\begin{align*}
\di \tilde \J^N(u)[h]&=\langle \di_x \J^N(S^N(u),u),\di S^N(u)[h]\rangle +\langle \di_u \J^N(S^N(u),u),h\rangle \\
&=\langle \di_u \J^N(S^N(u),u)-(\di_u e^N(S^N (u),u))^*(\di_x e^N(S^N (u),u))^{-*}\di_x \J^N(S^N(u),u),h\rangle .
\end{align*}
Defining the adjoint variable $\xi=(\xi_k)_{k=1}^N$ by
\begin{align*}
(\di_x e^N(S^N(u),u))^*[\xi]=-\di_x \J^N(S^N(u),u)
\end{align*}
yields
\begin{align*}
\di \tilde \J^N(u)[h]&=\langle \di_u \J^N(S^N(u),u)+(\di_u e^N(S^N (u),u))^*\xi,h\rangle\\
&= \dd \J_2(u)[h] +\dd \J_3(u)[h] 
\\&\quad- \frac{1}{N} \sum_{i=1}^N \sum_{j=1}^N \int_0^T 
\nabla_u \bar F(x_i(t) - x_j(t), u(x_i(t))) h(x_i(t))\cdot \xi_i(t)\di t.
\end{align*}
Using the variational lemma we can identify the gradient as
\begin{align}\label{eq:reducedcost}
\nabla \tilde J^N(u) = \dd \J_2(u) +\dd \J_3(u) 
- \frac{1}{N} \sum_{i=1}^N \sum_{j=1}^N \int_0^T 
\nabla_u \bar F(x_i(t) - x_j(t), u(x_i(t))) \xi_i(t)\di t.
\end{align}
With the gradient of the reduced cost functional \eqref{eq:reducedcost} at hand, we have everything required to state the gradient descent algorithm used for the computation of optimal controls.

\subsection{Convergence of the discrete optimal control problem}
As a first step towards the convergence of the discrete optimal control problem, we consider a stability estimate of the solutions to the discrete and the continuous adjoint problems \eqref{eq:adjointcont} and \eqref{eq:adjointdiscrete}. Similarly as in \cite{BPTT}, one can show the following stability estimate for the adjoint solution:
\begin{lmm}
	Let $x^N=(x_1,\ldots,x_N)$ be the solution to the forward particle system \eqref{eq:discretemodel} with initial condition \eqref{eq:discreteinitial} and given control $u^N\in\Uad$. Let $\rho^N(t,\cdot)$ denote the empirical measure corresponding to $x^N(t)$ for any $t\in [0,T]$. Let $\rho\in C([0, T ], \mathcal P_2(\Omega))$ be the solution to the mean-field state problem \eqref{eq:macroscopiceq} for given control $\bar u\in \Uad$. Let $\xi^N=(\xi_1,\ldots,\xi_N)$ denote the solution to the discrete adjoint system \eqref{eq:adjointdiscrete} for the pair $(x^N,u^N)$ and suppose that $\bar\xi\in C([0,T],\emph{Lip}_b(\Omega))$ satisfies \eqref{eq:adjointcont} for $(\rho,\bar u)$. Then, there exist positive constants $a$ and $b$, independent of $N\in \mathbb{N}$ such that
	\begin{align*}
		\sup_{t\in [0,T]} \frac{1}{N}\sum_{i=1}^N |\xi_i(t)-\bar\xi_t(x_i(t))|\leq b \exp(aT) \int_0^T  W_2(\rho^N(s,\cdot),\rho(s,\cdot)) +\|u^N-\bar u\|_\infty\di s.
	\end{align*} 
\end{lmm}	
Denoting by  $\rho^N(0,\cdot)$ the empirical measure which corresponds to the particle locations in $x^N(t)$ at time $t$ and using \eqref{eq:W2bound}, we obtain
\begin{align}\label{eq:convergencexi}
	\sup_{t\in [0,T]} \frac{1}{N}\sum_{i=1}^N |\xi_i(t)-\bar\xi_t(x_i(t))|^2\leq 
C_T\left(W_2^2(\rho^N(0,\cdot),\rho^\text{in}) + \|u-\bar{u}\|^2_\infty\right)
\end{align} 
for some constant $C_T$, depending on $T>0$ and independent of $N\in\mathbb{N}$.

\begin{thrm}
Let $(\bar \rho,\bar u)$ and $(x^N,u^N)$ be the optimal pairs for Problem \ref{prob:macroscopic} and Problem~\ref{prob:discrete} with  initial data $\rho^\text{in}\in \mathcal{P}(\Omega)$ and $X^N=(X_1,\ldots,X_N)\in \Omega^N$, respectively. Let $\rho^N(0,\cdot)$ denote the empirical measure corresponding to the initial configuration $X^N$. 
Then, there exists a constant $c_0 >0$ depending only on $T$,  $\bar \xi$ and $\bar F$   such that for $\lambda_1,\lambda_2>c_0$ in $\J_2,\J_3$, defined in \eqref{eq:deltaJ_2},  it holds 
\begin{align*}
\| u^N-\bar u\|_\infty^2\leq \frac{c_0}{\min\{\lambda_1,\lambda_2\}-c_0} W_2^2(\rho^N(0,\cdot),\rho^\text{in}).
\end{align*}
\end{thrm}	
\begin{proof}
	Let $\xi^N=(\xi_1,\ldots,\xi_N)$ denote the solution to the discrete adjoint system \eqref{eq:adjointdiscrete} for the pair $(x^N,u^N)$ and suppose that $\bar \xi\in C([0,T]\times \Omega)$ satisfies \eqref{eq:adjointcont} for $(\rho,u)$. We denote the empirical measure by $\rho_t^N$ which corresponds to the particle locations in $x^N(t)$ at time $t$. Considering the optimality conditions  \eqref{eq:optimality3} and \eqref{eq:optimalitydiscrete} in the macroscopic and microscopic setting, we have for $h^N:=u^N-\bar u\in U$
	 \begin{align*}
	  &(\dd \J_2(u^N)
	 - \dd \J_2(\bar u))[h^N]+(\dd \J_3(u^N)
	 - \dd \J_3(\bar u))[h^N] \\&=\frac{1}{N^2} \sum_{i=1}^N \sum_{j=1}^N \int_0^T 
	 \nabla_u \bF(x_i(t) - x_j(t), u^N(x_i(t))) h^N(x_i(t))\cdot \xi_i(t)\di t \\
	  &\quad- \int_0^T\int_\Omega \int_\Omega \nabla_u \bar F(x-y, \bar u(x))h^N(x) \cdot \bar \xi_t(x) \dd \bar \rho_t(y)\dd \bar \rho_t(x) \dd t\\&=\frac{1}{N}\sum_{i=1}^N\int_0^T \int_\Omega
	  \nabla_u \bF(x_i(t)- y, u^N(x_i(t))) h^N(x_i(t)) \dd  \rho^N(y) \cdot \xi_i(t) \di t \\
	  &\quad- \int_0^T\int_\Omega \int_\Omega \nabla_u \bar F(x-y, \bar u(x))h^N(x)\dd \bar \rho_t(y) \cdot \bar \xi_t(x) \dd \bar \rho_t(x) \dd t
	  \\&=\frac{1}{N}\sum_{i=1}^N\int_0^T \int_\Omega
	  \nabla_u \bF(x_i(t)- y, u^N(x_i(t))) h^N(x_i(t)) \dd  \rho_t^N(y) \cdot (\xi_i(t)-\bar \xi_t(x_i(t))) \di t\\
	  &\quad+\int_0^T \int_\Omega\int_\Omega
	  \nabla_u \bF(x- y, u^N(x)) h^N(x)\cdot \bar \xi_t(x) \dd  \rho_t^N(x)\dd  \rho_t^N(y)\di t \\
	  &\quad-\int_0^T \int_\Omega\int_\Omega\nabla_u \bar F(x-y, \bar u(x))h^N(x) \cdot \bar \xi_t(x)\dd \bar \rho_t(x) \dd \bar \rho_t(y)\di t.
	  \end{align*} 
	 Using \eqref{eq:convergencexi}, the first term can be estimated by 
	 \begin{align*}
	  C_T\left(W_2^2(\rho^N(0,\cdot),\rho^\text{in}) + \|u^N-\bar{u}\|^2_\infty\right)
	 \end{align*}
	 for some constant $C_T$ depending of $T$ and independent of $N$. Denoting the optimal coupling between $\rho_t^N$ and $\bar \rho_t$ by $\pi_t$, the remaining terms can be rewritten as
	 	  \begin{align*}
	 	  & \int_0^T\int_\Omega \int_{\Omega\times \Omega} \left( \nabla_u \bF(x- y, u^N(x)) h^N(x)\cdot \bar \xi_t(x) - \nabla_u \bar F(x'-y,  u^N(x'))h^N(x') \cdot \bar \xi_t(x') \right)\dd  \pi_t(x,x') \dd  \rho^N_t(y)\dd t\\
	 	  &\quad+ \int_0^T\int_\Omega \int_{\Omega} \left[ \nabla_u \bF(x- y, u^N(x))  - \nabla_u \bar F(x-y, \bar u(x))\right]h^N(x) \cdot \bar \xi_t(x)  \dd \bar \rho_t(x) \dd  \rho^N_t(y)\di t\\
	 	  &\quad+ \int_0^T\int_{\Omega\times \Omega} \int_{\Omega} \left[ \nabla_u \bF(x- y, \bar u(x)) - \nabla_u \bar F(x-y', \bar u(x)) \right] h^N(x) \cdot \bar \xi_t(x)  \dd \bar  \rho_t(x) \dd  \pi_t(y,y')\dd t\\
	 	  &\leq c_T \left( \left(	\|\nabla_u \bar  F\|_{\infty} \sup_{t\in[0,T]} \operatorname{Lip}(\bar\xi_t)+\|\nabla_x\nabla_u \bar  F\|_{\infty} \sup_{t\in[0,T]} \| \bar\xi_t\|_{\infty} \right)\|u^N-\bar{u}\|_\infty  W_2(\rho^N_t,\bar\rho_t)\right.\\
	 	  &\quad+\left.\|\nabla_u \nabla_u \bar  F\|_{\infty} \sup_{t\in[0,T]} \| \bar\xi_t\|_{\infty} \|u^N-\bar{u}\|_\infty^2  +\|\nabla_x\nabla_u \bar  F\|_{\infty} \sup_{t\in[0,T]} \| \bar\xi_t\|_{\infty} \|u^N-\bar{u}\|_\infty  W_2(\rho^N_t,\bar\rho_t)\right)
	 \end{align*}
	for some constant $c_T$. This yields 
	\begin{align*}
	(\dd \J_2(u^N)
	- \dd \J_2(\bar u))[h^N] +	(\dd \J_3(u^N)
	- \dd \J_3(\bar u))[h^N]\leq c_0\left(	\|u^N-\bar{u}\|^2_\infty + W_2^2(\rho^N_t,\bar\rho_t)\right)
	\end{align*}
	where $c_0$ depends on $T$, $\bar\xi$ and $\bar F$. For $u^N=(\theta^N,\eta^N)$, $\bar u=(\bar \theta,\bar \eta)$ and $h^N=(\theta^N-\bar \theta,\eta^N-\bar \eta)$, we have 
	\begin{align*}
	(\dd \J_2(u^N)
	- \dd \J_2(\bar u))[h^N]+(\dd \J_3(u^N)
	- \dd \J_3(\bar u))[h^N] & = \lambda_1  \| \theta^N-\bar \theta\|_{\mathcal H(\Omega)}^2  + \lambda_2  \| \eta- \bar \eta\|_{\mathcal H(\Omega)}^2\\
	&\geq \min\{\lambda_1,\lambda_2\} \|u^N-\bar{u}\|^2_\infty
	\end{align*}
	by \eqref{eq:deltaJ_2} and the continuous embedding of $\mathcal H(\Omega)$ in $L^\infty(\Omega)$.
	This implies
	\begin{align*}
	( \min\{\lambda_1,\lambda_2\}-c_0) \|u^N-\bar{u}\|^2_\infty \leq  c_0 W_2^2(\rho^N_t,\bar\rho_t).
	\end{align*}	
	Choosing $\lambda_1,\lambda_2>c_0$ we obtain the desired inequality.
	\end{proof}

\section{Numerical schemes}\label{sec:numericalschemes}
In this section, we introduce the numerical schemes, used for solving the forward and adjoint initial value problems including the terminal condition of the adjoint problem, as well as the optimal control problem.

\subsection{Forward and adjoint initial value problems}
For solving the discrete forward problem \eqref{eq:discretemodel} with initial condition \eqref{eq:discreteinitial} on the unit square $[0,1]^2$ we apply the simple explicit Euler scheme. The discrete adjoint problem \eqref{eq:adjointdiscrete} is linear and stiff, it is therefore solved implicitly. Note that the force $\bF$ is defined periodically by \eqref{eq:forceperiodic} in both problems. Further, periodic boundary conditions guarantee that the particle positions $x_i$ cannot leave the domain, i.e.\ $x_i\in[0,1]^2$ for $i=1,\ldots,N$.

\begin{rmrk}
We emphasize that the first-optimize then discretize approach allows us to choose different discretizations for the forward and the adjoint solver. This is a huge advantage, as otherwise the computational effort increases tremendously due to very small step sizes or a complicated implementation of the forward solver.
\end{rmrk}

\subsection{Terminal condition of the adjoint problem}
The main challenge of the implementation of the particle optimization is the evaluation of the terminal condition of the adjoints, given by
\begin{align}\label{eq:terminalcondadjoint}
\xi_k(T):= \bar \xi_T(x_k) = N \delta_\rho \J_1(\bar \rho_T)(x_k) =N  t_{\rho_T}^{\rho_\text{des}} (x_k(T)) -x_k(T). 
\end{align} 
We realize it with the help of the \texttt{Python Optimal transport} library \cite{pythonOT}.

While $\rho_T$ and $\rho_\text{des}$ are probability densities in the macroscopic setting, we consider the associated empirical measures $\rho^N(T,\cdot)=\frac{1}{N}\sum_{k=1}^N \delta_{x_k(T)}$ and $\rho_\text{des}^N=\frac{1}{N}\sum_{k=1}^N \delta_{x_k^\text{des}}$, respectively. Let $\texttt{a,b}$ denote the sample weights for the 1D histograms corresponding to $x_k(T)$ and $x_k^\text{des}$ for $k=1,\dots,N,$ i.e.\ $\texttt{a}_k=\frac{1}{N}=\texttt{b}_k$ for $k=1,\ldots,N$. Instead of the usual ground cost matrix, we use a ground matrix $\texttt{M}$ that accounts for the periodic boundary conditions. 
We compute the earth mover's distance (EMD) using the function \texttt{G0 = ot.emd(a,b,M)} where $\texttt{G0}$ solves  Kantorovich's optimal transport problem:
\begin{align*}
\texttt{G0 } &=\argmin_{\texttt{G0}\in \texttt{U(a,b)}}\langle\texttt{G0,M} \rangle
\end{align*}
for $$\texttt{U(a,b)}=\{ \texttt{G}\in \R^{N,N} \colon \texttt{G1\textsubscript{N}=a, G\textsuperscript{T}1\textsubscript{N}=b} \},$$
where $\texttt{1\textsubscript{N}}$ is the vector of ones of length $N$.
In particular, each entry $\texttt{G0\textsubscript{ij}}$ of the coupling matrix $\texttt{G0}$ describes the amount of mass flowing from the mass found at $x_i(T)$ towards $x^{\text{des}}_{j}$.
  
  Having $\texttt{G0}$ at hand, we compute the vectors connecting the $x_k(T)$ with the desired position $t_{\rho_T}^{\rho_\text{des}}(x_k(T))$ for all $k=1,\dots,N.$ To do this, we assemble position vectors $\vec x,\vec y\in \R^N$ containing the first and the second component of the positions of the particles, respectively. We proceed analogously with the positions of the particles of the desired distribution. Then we construct the matrices
\begin{align*} X_1 &= (\vec x, \vec x, \dots, \vec x), \quad X_2 = (\vec x_\text{des}, \vec x_\text{des}, \dots, \vec x_\text{des})^T, \\  Y_1 &= (\vec y, \vec y, \dots, \vec y),\quad Y_2 = (\vec y_\text{des}, \vec y_\text{des}, \dots, \vec y_\text{des})^T.
\end{align*}
The distances between each particle of the current and the desired distribution, again accounting for the periodic boundary conditions, are contained in the matrices
$$ M_X = X_2-X_1, \qquad M_Y = Y_2-Y_1. $$
For the right-hand side of the adjoints we use
\[
r_x = \text{sum}(M_X*(\texttt{G0}>0)), \qquad r_y = \text{sum}(M_Y*(\texttt{G0}>0)),
\]
where $*$ denotes componentwise multiplication and entries of $\texttt{G0}$ are only considered for $\texttt{G0}>0$.

\subsection{Optimal control problem}
As mentioned before, the simulation results are very sensitive to the value of $\eta.$ We therefore restrict the domain to $\eta \in [\eta_\text{min}, \eta_\text{max}].$ Thus, given a control $u^N$, we update the control by $v^N$ via a projected gradient decent \cite{Hinze}, where the  step size $\tau$ is determined via line search, i.e.\ we consider
\begin{align*}
v^N= \mathcal{P}_{\Uad} \big(u^N+\tau \nabla \tilde J^N(u^N) \big),
\end{align*}
where $\nabla \tilde J^N(u^N)$ denotes the gradient of the reduced cost functional in \eqref{eq:reducedcost} and $\mathcal{P}_{\Uad}$ is the projection onto $\Uad$.
Using the solvers for the state system \eqref{eq:discretemodel}, the adjoint system \eqref{eq:adjointdiscrete} and steepest descent to update the control $u^N$, we obtain Algorithm \ref{alg:controlalg}.

 \begin{algorithm}[H]
	\SetAlgoLined
	\KwData{Initial data $x^N(0)=(x_1(0),\ldots,x_N(0))$ with $x_i(0)\in[0,1]^2$ given by \eqref{eq:discreteinitial} for $i=1,\ldots,N$; simulation time $T>0$; desired values $x^N_\text{des}$; other parameter values;}
	\KwResult{Control $u^N$,  state $x^N(T)$ and optimal function value $\J^N(x^N(T),u^N)$}
	initialization\;
	solve state problem \eqref{eq:discretemodel} for $x^N$\;
	solve adjoint problem \eqref{eq:adjointdiscrete} for $\xi^N$\;
	evaluate gradient of the reduced cost functional $\nabla \tilde J^N(u^N)$ in \eqref{eq:reducedcost}\;
	\While{stopping criterion not satisfied}{
		perform a line search to update the control $u^N$, compute $x^N$ and evaluate $\J^N$\;
		solve adjoint problem \eqref{eq:adjointdiscrete} for $\xi^N$\;
		evaluate gradient of the reduced cost functional $\nabla \tilde J^N(u^N)$ in \eqref{eq:reducedcost}\;
	}
	\caption{Optimal Control Algorithm}\label{alg:controlalg}
\end{algorithm}

\section{Numerical results}\label{sec:numericalresults}
In this section we discuss numerical results, obtained with the particle algorithm introduced in Algorithm~\ref{alg:controlalg}. Since we are mainly interested in the Wasserstein distance in the cost functional, we restrict ourselves to the simple case of spatially homogeneous control parameters, i.e., $u=(\theta,\eta) \in \R \times \R$, for the numerical simulations. Note that the norms in $J_2$ and $J_3$ reduce to the standard norms in $\R.$ 

First, we set the force coefficients and the parameter values. Then we describe how  artificial data is obtained for the parameter estimation, as we have no real data available. Finally, we show simulation and convergence results.
\subsection{Parameter values for the results}
For the numerical examples we set $N=1200$ and choose the force coefficients $f_R$ and $f_A$ of the repulsion and attraction forces \eqref{eq:repulsionforce} and \eqref{eq:attractionforce} as
\begin{equation*}
f_R(\eta |d|) = (\alpha \eta^2 |d|^2 + \beta) e^{-e_R \eta |d|}, \qquad f_A(\eta |d|) = -\gamma \eta |d| e^{- e_A \eta|d|}, 
\end{equation*}
resulting in the total forces
\begin{equation*}
F_R(d=d(x,y),u)=\eta f_R(\eta |d|)d, \qquad
F_A(d=d(x,y),u)=\eta f_A(\eta |d|) R_\theta \begin{pmatrix} 1 & 0 \\ 0 & \chi\end{pmatrix} R_\theta^T d,
\end{equation*}
with parameters $\alpha = 270, \beta = 0.1, \gamma= 35, e_R = 100, e_A = 95$ as in \cite{patternformationanisotropicmodel}. 
Moreover, we set $dt = 2$ and $T=10000$ leading to $5000$ time steps for each solve of the forward, adjoint and gradient computation. The optimization iteration is stopped when the relative gradient satisfies the condition
\[
\frac{\| \nabla J_k \|_1}{\| \nabla J_0 \|_1} < \epsilon_\text{stop},
\]
where $\nabla J_0$ corresponds to the first gradient of the computation and $\nabla J_k$ denotes the gradient of the current iteration in the optimization procedure. We choose $\epsilon_\text{stop} = 0.05$ for all simulations. The reference values in the cost function are set to $\theta_\text{ref} = 0.5\pi$ and $\eta_\text{ref} = 1$, and the scaling factors in the cost function are $\lambda_1 = 1e^{-5}$ and $\lambda_2=1e^{-3}.$  The parameter value for the anisotropy is $\chi = 0.2.$ The admissible interval for $\eta$ is given by $[\eta_\text{min}, \eta_\text{max}] = [0.9, 1.1].$

\subsection{Artificial data}
As we do not have any real data available, we compute some artificial data to validate our approach. We therefore choose the parameters $\theta_\text{data}$ and $\eta_\text{data}$, and consider some random initial condition $x(0)$ which consists of uniformly distributed positions in $[0,1]^2.$ The initial positions for the optimization is another sample of uniformly distributed positions in $[0,1]^2.$ This induces some noise and therefore, we do not expect that the algorithm fits the data perfectly. We use different values for the data parameters and give  details for every simulation together with the corresponding simulation results.

\subsection{Simulation results}
\begin{figure}[htp]
\hspace{-0.4cm}	\includegraphics[scale=0.8]{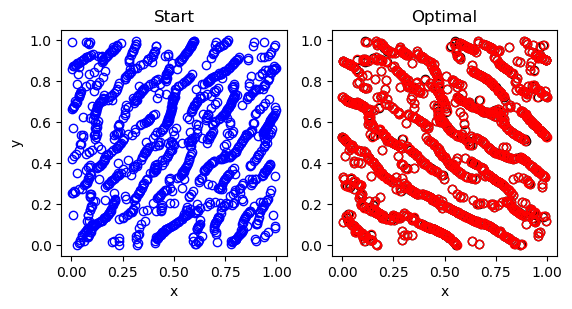}\\
	\includegraphics[scale=0.8]{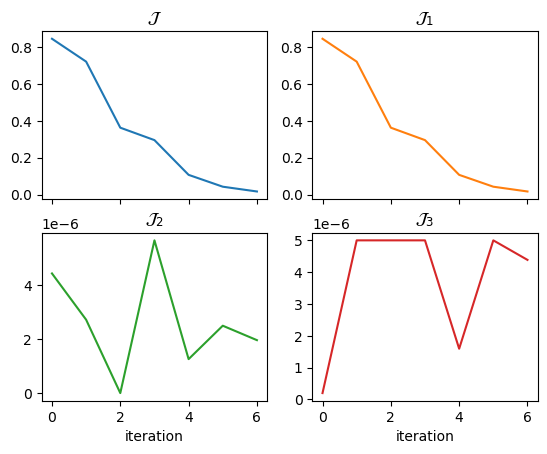}	
	\caption{P2) On the top-left we see the state at $T=10000$ for the initial values $\theta_0 = 0.8\pi$ and $\eta_0 = 0.98.$ The plot on the top-right reports the states of the artificial data at $T=10000$ in black and the state corresponding to the optimized values $\theta_\text{opt} =0.3003\pi, \eta_\text{opt}=0.9063$ in red. On the bottom, the evolution of the cost functional values for each optimization iteration is shown. }
	\label{fig:P2}
\end{figure}

\begin{figure}[htp]
\hspace{-0.4cm}	\includegraphics[scale=0.8]{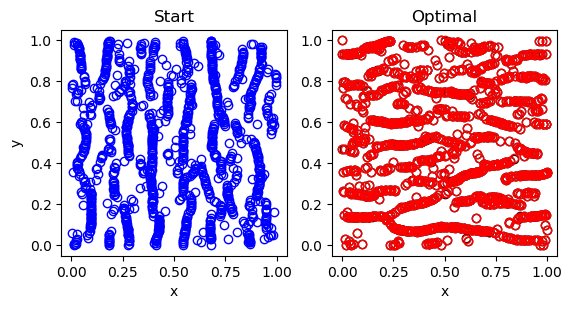}
	\includegraphics[scale=0.8]{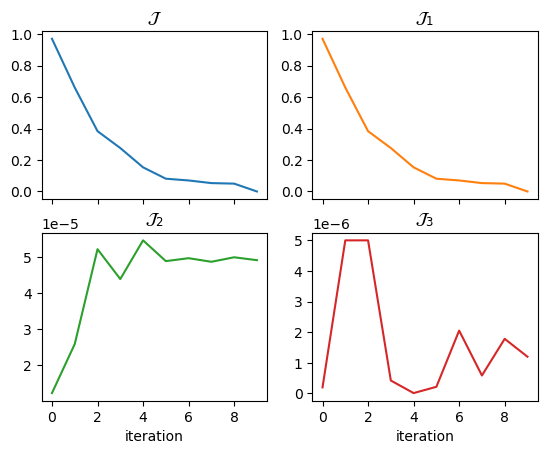}	
	\caption{P3) On the top-left we see the state at $T=10000$ for the initial values $\theta_0 = 0.0\pi$ and $\eta_0 = 0.98.$ The plot on the top-right reports the states of the artificial data at $T=10000$ in black and the state corresponding to the optimized values $\theta_\text{opt} = -0.4989\pi, \eta_\text{opt}=0.95099$ in red. On the bottom, the evolution of the cost functional values for each optimization iteration is shown. }
	\label{fig:P3}
\end{figure}

Let $\theta_0, \eta_0$ denote the initialization values for the optimization procedure  and let $\theta_\text{data}, \eta_\text{data}$ be the values for the artificial data. The simulation results shown below correspond to the following parameters:
\begin{itemize}
	\item [P1)] data: $\theta_0 = 0.3\pi,\; \theta_\text{data} = 0.7\pi,\; \eta_0 = 0.98, \;\eta_\text{data} = 1.0$, \\ optimized values: $\theta_\text{opt} = 0.7035\pi,\; \eta_\text{opt} = 1.0221$ 
	\item [P2)] data: $\theta_0 = 0.8\pi,\; \theta_\text{data} = 0.3\pi,\; \eta_0 = 0.98, \;\eta_\text{data} = 0.9$, \\ optimized values: $\theta_\text{opt} = 0.3003\pi,\, \eta_\text{opt} = 0.9063$ 
	\item [P3)] data: $\theta_0 = 0.0\pi,\; \theta_\text{data} = 0.5\pi,\; \eta_0 = 0.98,\;\eta_\text{data} = 0.95$, \\ optimized values: $\theta_\text{opt} = -0.4989\pi,\; \eta_\text{opt} = 0.95099$
\end{itemize}
\begin{figure}[htp]
	\hspace{-0.3cm}\includegraphics[scale=0.8]{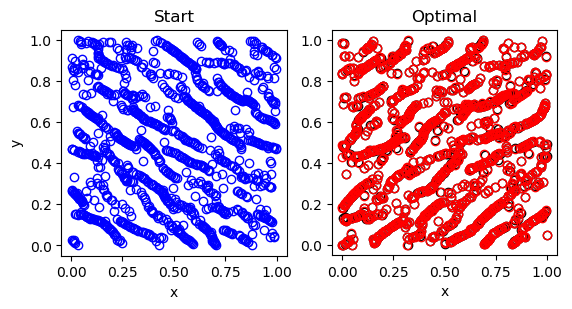}\\
	\includegraphics[scale=0.8]{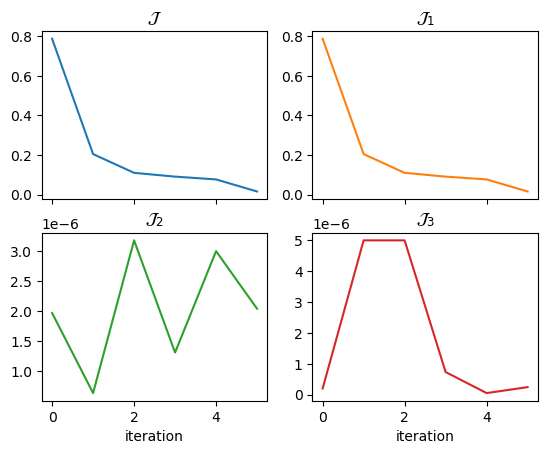}	
	\caption{P1) On the top-left we see the state at $T=10000$ for the initial values $\theta_0 = 0.3\pi$ and $\eta_0 = 0.98.$ The plot on the top-right reports the states of the artificial data at $T=10000$ in black and the state corresponding to the optimized values $\theta_\text{opt} =0.7035\pi, \eta_\text{opt}=1.0221$ in red. On the bottom, the evolution of the cost functional values for each optimization iteration is shown.}
	\label{fig:P1}
\end{figure}
The initial states and the ones corresponding to the optimized parameter values P1), P2), P3) are shown in Figure~\ref{fig:P1}, Figure~\ref{fig:P2} and Figure~\ref{fig:P3}, respectively. The results indicate a good performance of the algorithm. Indeed, in the eye-norm there is no difference in the states obtained with the artificial data and the states corresponding to the optimized parameters visible, as shown in the plots on the top-right of the Figures~\ref{fig:P1}-\ref{fig:P3}. Moreover, the maximum number of optimization iterations is $9$. The cost functional values decrease as usual in optimal control, that means, the first optimization steps reduce the cost more than later steps. As the regularization values $\lambda_1$ and $\lambda_2$ are chosen very small, the total cost is mainly driven by $\J_1$ which corresponds to the Wasserstein distance of the discrete densities. This is the desired behaviour. 

It is interesting to see that for $P3)$ the optimal angle is approximating the shifted reference angle, i.e., $\theta_\text{opt} \approx \theta_\text{art} - \pi$. This occurs  as we allow $\theta \in \R$ and have very small regularization parameters $\lambda_1$ and $\lambda_2.$ We also see the difference  in the plots of the cost functional $\mathcal J_2$. Indeed, note that for $P3)$ the value of $\mathcal J_2$ at the end of the optimization is about one order of magnitude larger then the ones corresponding to $P1)$ and $P2).$

\section{Conclusion}\label{sec:conclusion}
We proposed a mean-field optimal control ansatz to identify parameters underlying given, artificially generated patterns. The state system is an agent-based model with anisotropic interaction forces that lives on the torus. The identification algorithm used gradient information that is computed with the help of the first order optimality conditions. The cost functional penalizes the Wasserstein distance of the data pattern and the modelled pattern resulting from the state system for large times. Numerical results on the particle level demonstrate the performance of the proposed method. These results can be seen as a first step towards the modelling of complex fingerprint patterns with specific features in future work.

\section*{acknowlegements} \noindent
MB has been partially supported by the German Science Foundation (DFG) through CRC TR 154 ”Mathematical Modelling, Simulation and Optimization Using the Example of Gas Networks”. MB and LMK acknowledge support from the European Union Horizon 2020 research and innovation programmes under the Marie Sk\l odowska-Curie grant agreement No.\ 777826 (NoMADS). LMK acknowledges support from the European Union Horizon 2020 research and innovation programmes under the Marie Sk\l odowska-Curie grant agreement No.\ 691070 (CHiPS), the EPSRC grant EP/L016516/1, the German National Academic Foundation (Studienstiftung des Deutschen Volkes), the Cantab Capital Institute for the Mathematics of Information and Magdalene College, Cambridge (Nevile Research Fellowship).
CT was partly supported by the European Social Fund and by the Ministry Of Science, Research and the Arts Baden-W\"urttemberg. Moreover, CT acknowledges support by the state of Baden-Württemberg through bwHPC, in particular the bwForCluster MLS\&WISO Production. 
	
\bibliographystyle{plain}
\bibliography{references}
	
\end{document}